\documentclass[12pt]{amsart}
\usepackage{amssymb}
\newtheorem{theorem}{Theorem}
\newtheorem{proposition}[theorem]{Proposition}

\def\ds{\displaystyle}

\def\Ree{\operatorname{Re}}
\def\Imm{\operatorname{Im}}
\title{On the derivatives of the Lempert functions}

\author{Nikolai Nikolov and Peter Pflug}

\address
{Institute of Mathematics and Informatics\\ Bulgarian Academy of
Sciences\\ Acad. G. Bonchev 8, 1113 Sofia,
Bulgaria}\email{nik@math.bas.bg}

\address{Carl von Ossietzky Universit\"at Oldenburg\\
Institut f\"ur Mathematik, Fakult\"at V\\ Postfach 2503\\ D-26111
Oldenburg, Germany}\email{pflug@mathematik.uni-oldenburg.de}



\begin{document}

\begin{thanks}{This note was written during the stay of the first
named author at the Universit\"at Oldenburg supported by a grant
from the DFG, Az. PF 227/8-2 (November -- December 2006). He likes
to thank both institutions for their support.}
\end{thanks}

\begin{abstract} We show that if the Kobayashi--Royden metric of a
complex manifold is continuous and positive at a given point and any
non-zero tangent vector, then the "derivatives" of the higher order
Lempert functions exist and equal the respective Kobayashi metrics
at the point. It is a generalization of a result by M.~Kobaya-\\shi
for taut manifolds.
\end{abstract}

\maketitle

\section{Introduction and results}

Let $\Bbb D\subset\Bbb C$ be the unit disc. Let $M$ be an
$n$-dimensional complex manifold. Recall first the definitions of
the Lempert function $\tilde k_M$ and the Kobayashi--Royden
pseudometric $\kappa_M$ of $M$:
$$\aligned
\tilde k^\ast_M(z,w)&=\inf\{|\alpha|:\exists f\in\mathcal
O(\Bbb D,M):f(0)=z,f(\alpha)=w\},\\
\tilde k_M&=\tanh^{-1}\tilde k^\ast_M,\\
\kappa_M(z;X)&=\inf\{|\alpha|:\exists f\in\mathcal O(\Bbb D,M):
f(0)=z,\alpha f_\ast(d/d\zeta)=X\},
\endaligned$$
where $X$ is a complex tangent vector to $M$ at $z$. Note that such
an $f$ always exists (cf.~\cite{Win}; according to \cite{Din}, page
49, this was already known by J. Globevnik).

The Kobayashi pseudodistance $k_M$ can be defined as the largest
pseudodistance bounded by $\tilde k_M$. Note that if $k_M^{(m)}$
denotes the $m$-th Lempert function of $M$, $m\in\Bbb N$, that is,
$$k_M^{(m)}(z,w)=\inf\{\sum_{j=1}^m\tilde k_M(z_{j-1},z_{j}):z_0,\dots,z_{m}\in M,
z_0=z,z_m=w\},$$ then
$$k_M(z,w)=k_M^{(\infty)}:=\inf_m k_M^{(m)}(z,w).$$

By a result of M.-Y.~Pang (see \cite{Pang}), the Kobayashi--Royden
metric is the "derivative" of the Lempert function for taut domains
in $\Bbb C^n;$ more precisely, if $D\subset\Bbb C^n$ is a taut
domain, then
$$\kappa_D(z;X)=\lim_{t\nrightarrow 0}\frac{\tilde k_D(z,z+tX)}{t}.$$

In \cite{KobS}, S.~Kobayashi introduces a new invariant
pseudometric, called the Kobayashi--Buseman pseudometric in
\cite{Jar-Pfl}. One of the equivalent ways to define the
Kobayashi--Buseman pseudometric $\hat\kappa_M$ of $M$ is just to set
$\hat\kappa_M(z;\cdot)$ to be largest pseudonorm  bounded by
$\kappa_M(z;\cdot)$. Recall that
$$
\hat\kappa_M(z;X)=\inf\{\sum_{j=1}^m\kappa_M(z;X_j):m\in\Bbb N,\
\sum_{j=1}^mX_j=X\}.
$$ Thus it is natural to consider the new function
$\kappa_M^{(m)}$, $m\in\Bbb N$, namely,
$$
\kappa_M^{(m)}(z;X)=\inf\{\sum_{j=1}^m\kappa_M(z;X_j):
\sum_{j=1}^mX_j=X\}.
$$
We call $\kappa_M^{(m)}$ the {\it $m$-th Kobayashi pseudometric} of
$D$. It is clear that $\kappa_M^{(m)}\ge\kappa_M^{(m+1)}$ and if
$\kappa_M^{(m)}(z;\cdot)=\kappa_M^{(m+1)}(z;\cdot)$ for some $m$,
then $\kappa_M^{(m)}(z;\cdot)\\=\kappa_D^{(j)}(z;\cdot)$ for any
$j>m$. It is shown in \cite{Nik-Pfl} that
$\kappa_M^{(2n-1)}=\kappa_M^{(\infty)}:=\hat\kappa_M$, and $2n-1$ is
the optimal number, in general.

We point out that all the introduced objects are upper
semicontinuous. Recall that this is true for $\kappa_M$ (cf.
\cite{KobS2}). It remains to check this for $\tilde k_M$. We shall
use a standard reasoning. Fix $r\in(0,1)$ and $z,w\in M$. Let
$f\in\mathcal O(\Bbb D,M)$, $f(0)=z$ and $f(\alpha)=w$. Then $\tilde
f=(f,\hbox{id}):\Delta\to\tilde M=M\times\Delta$ is an embedding.
Setting $\tilde f_r(\zeta)=\tilde f(r\zeta)$, by \cite{Roy}, Lemma
3, we may find a Stein neighborhood $S\subset\tilde M$ of $\tilde
f_r(\Bbb D)$. Embed $S$ as a closed complex manifold in some $\Bbb
C^N$ and denote by $\psi$ the respective embedding. Moreover, there
is an open neighborhood $V\subset \Bbb C^N$ of $\psi(S)$ and a
holomorphic retraction $\theta:V\to\psi(S)$. Then, for $z'$ near $z$
and $w'$ near $w$, we may find, as usual, $g\in\mathcal O(\Bbb D,V)$
such that $g(0)=\psi(z',0)$ and $g(\alpha/r)=\psi(w',\alpha)$.
Denote by $\pi$ the natural projection of $\tilde M$ onto M. Then
$h=\pi\circ\psi^{-1}\circ\theta\circ g\in\mathcal O(\Bbb D,M)$,
$h(0)=z'$ and $h(\alpha/r)=w'$. So $r\tilde
k^\ast_M(z',w')\le\alpha$, which implies that $\ds\limsup_{z'\to
z,w'\to w}\tilde k_M(z',w')\le\tilde k_M(z,w)$.

To extend Pang's result on manifolds, we have to define the
"derivati\-ves" of $k^{(m)}_M$, $m\in\Bbb N^\ast=\Bbb
N\cup\{\infty\}$. Let $(U,\varphi)$ be a holomorphic chart near
$z$. Set $$\mathcal D k^{(m)}_M(z;X)=\limsup_{t\nrightarrow 0,w\to
z,Y\to\varphi_\ast X}
\frac{k^{(m)}_M(w,\varphi^{-1}(\varphi(w)+tY))}{|t|}.$$ Note that
this notion does not depend on the chart used in the definition
and $$\mathcal D k^{(m)}_M(z;\lambda X)=|\lambda|\mathcal D
k^{(m)}_M(z;X),\quad \lambda\in\Bbb C.$$ Replacing $\limsup$ by
$\liminf$, we define $\underline{\mathcal D} k_M^{(m)}$.

From M. Kobayashi's paper \cite{KobM} it follows that, if $M$ is a
taut manifold, then
$$
\hat\kappa_M(z;X)=\mathcal D k_M(z;X)=\underline{\mathcal D}
k_M(z;X),$$ that is, the Kobayashi--Buseman metric is the
"derivative" of the Koba\-yashi distance. The proof there also leads
to $$\kappa_M^{(m)}(z;X)=\mathcal D
k^{(m)}_M(z;X)=\underline{\mathcal D} k_M^{(m)}(z;X),\quad m\in\Bbb
N^{\ast}.\leqno{(\ast)}$$

We say that a complex manifold $M$ is hyperbolic at $z$ if
$k_M(z,w)>0$ for any $w\neq z$. We point out that the following
conditions are equivalent:

(i) M is hyperbolic at $z;$

(ii) $\ds\liminf_{z'\to z,w\in M\setminus U}\tilde k_M(z',w)>0$ for
any neighborhood $U$ of $z;$

(iii) $\ds\underline{\kappa}_M(z;X):=\liminf_{z'\to z,X'\to
X}\kappa_M(z';X')>0$ for any $X\neq 0;$

The implication (i)$\Rightarrow$(ii) $\Rightarrow$(iii) are almost
trivial (cf. \cite{Jar-Pfl}) and the implication
(iii)$\Rightarrow$(i) is a consequence of the fact that $k_M$ is the
integra\-ted form of $\kappa_M$.

In particular, if $M$ is hyperbolic at $z$, then it is hyperbolic at
any $z'$ near $z$.

Since if $M$ is taut, then it is k-hyperbolic and $\kappa_M$ is a
continuous function, the following theorem is a generalization of
$(\ast)$.

\begin{theorem}\label{th1} Let $M$ be a complex manifold and $z\in M$.

{\rm (i)} If $M$ is hyperbolic at $z$ and $\kappa_M$ is continuous
at $(z,X)$, then
$$\kappa_M(z;X)=\mathcal D\tilde k_M(z;X)=\underline{\mathcal D}
\tilde k_M(z;X).$$

{\rm(ii)} If $\kappa_M$ is continuous and positive at $(z,X)$ for
any $X\neq 0$, then
$$\kappa_M^{(m)}(z;\cdot)=\mathcal D k_M^{(m)}(z;\cdot)=\underline{\mathcal D}
k_M^{(m)}(z;\cdot),\quad m\in\Bbb N^\ast.$$

\end{theorem}

The first step in the proof of  Theorem \ref{th1} is  the following

\begin{proposition}\label{pr2} For any complex manifold $M$ one has
that
$$\kappa_M^{(m)}\ge\mathcal D k_M^{(m)}, \quad m\in\Bbb N^\ast.$$
\end{proposition}

Note that when $M$ is a domain, a weaker version of Proposition
\ref{pr2} can be found in \cite{Jar-Pfl}, namely,
$\hat\kappa_M\ge\mathcal D k_M$ (the proof is based on the fact that
$\mathcal D k_M(z;\cdot)$ is a pseudonorm).
\section{Examples}

The following examples show that the assumption on continuity in
Theorem \ref{th1} is essential.

\smallskip

\noindent$\bullet$ Let $A$ be a countable dense subset of $\Bbb
C_\ast.$ In \cite{Die-Sib} (see also \cite{Jar-Pfl}), a pseudoconvex
domain $D$ in $\Bbb C^2$ is constructed such that:

(a) $(\Bbb C\times\{0\})\cup(A\times\Bbb C)\subset D;$

(b) if $z_0=(0,t)\in D,$ $t\neq 0,$ then $\kappa_D(z_0;X)\ge C||X||$
for some $C=C_t>0$. (One can be shown that even $\mathcal D\tilde
k_D(z_0;X)\ge C||X||.$)

\noindent Then it is easy to see that $\underline{
\kappa}_D(\cdot;e_2)=\mathcal D k^{(3)}_D(\cdot;e_2)=k^{(5)}_D=0$
and $\hat\kappa_D(z_0;X)\ge c||X||,$ where $e_2=(0,1)$ and $c>0.$
Thus $$\hat\kappa_D(z_0;X)>\underline{ \kappa}_D(z_0;e_2)=\mathcal
D k^{(3)}_D(z_0;e_2)=\mathcal D k^{(5)}_D(z_0;X),\quad X\neq 0.$$

This phenomena obviously extends to $\Bbb C^n$, $n>2$ (by
considering $D\times\Bbb D^{n-2}$). So the inequalities in
Proposition \ref{pr2} are strict in general.

\smallskip

\noindent$\bullet$ If $D$ is a pseudoconvex balanced domain with
Minkowski function $h_D$, then (cf. \cite{Jar-Pfl})
$$h_D=\kappa_D(0;\cdot)=\mathcal D{\tilde k}_D(0;\cdot).$$
Therefore, $\mathcal D{\tilde k}_D(0;X)>\underline{\mathcal
D}\tilde k_D(0;X)$ if $\kappa_D(0;\cdot)$ is not continuous at
$X$. On the other hand, if $\hat D$ denotes the convex hull of
$D$, then $$h_{\hat D}=\hat\kappa_D(0;\cdot)=\mathcal
D{k}_D(0;\cdot)=\underline{\mathcal D}k_D(0;\cdot)=
\underline{\hat\kappa}_D(0;\cdot).$$

\smallskip

\noindent$\bullet$ Modifying the first example leads to a
pseudoconvex domain $D\subset\Bbb C^2$ with $$L_{\mathcal
Dk_D}(\gamma)>0=L_{k_D}(\gamma)= L_{\underline{\mathcal
D}k_D}(\gamma),$$ where $\gamma:[0,1]\to\Bbb C^2$,
$\gamma(t):=(ti/2,1/2)$, and $L_\bullet(\gamma)$ denotes the
respective length.

Indeed, choose a dense sequence $(r_j)$ in $[0,i/2].$ Put $$
u(\lambda)=\sum_{k=1}^\infty
\frac{1}{k^2}\log\tfrac{|\lambda-1/k|}{4},\quad
v(\lambda)=\sum_{j=1}^\infty\frac{u(\lambda/2-r_j)}{2j^2},\quad\lambda\in\Bbb
C,
$$
and
$$
 D=\{z\in\Bbb C^2:\psi(z)=|z_2|e^{\|z\|^2+v(z_1)}<1\}.
$$
It is easy to see that $v$ is a subharmonic function on $\Bbb C$.
Hence $D$ is a pseudoconvex domain with $(\Bbb C
\times\{0\})\cup(\bigcup_{j,k=1}^\infty \{r_j+1/k\}\times\Bbb
C)\subset D.$ Observe that $u|_{\Bbb D}<-1$ and so $D$ contains the
unit ball $\Bbb B_2.$ Note also that
$$
k_D(a,b)=0,\quad a,b\in\gamma([0,1]).
$$
Set $\widehat \psi(z)=\|z\|^2/2-\log\psi(z)$. Fix $z^0\in\Bbb B_2$
with $\Ree z^0_1\le 0,\;\Imm z^0_2\ge1/e$. Since $u(\lambda)\ge
u(0)$ for $\Ree \lambda\le 0,$ we have
$$
||z^0||/2<\widehat\psi(z^0)<1-u(0)=:8C.
$$

Let $\varphi\in\mathcal O(\Bbb D,D)$, $\varphi(0)=z^0.$ Following
the estimates in the proof of Example 3.5.10 in \cite{Jar-Pfl}, we
see that $\|\varphi'(0)\|<C$. Hence, $\kappa_D(z^0;X)\geq C\|X\|$,
$X\in\Bbb C^2$. Since $k_D$ is the integrated form of $\kappa_D,$
it follows that $$ k_D(a,a-te_1)\geq Ct,\quad
a\in\gamma([0,1]),\;0\le t\le 1/2-1/e,\;e_1=(1,0).$$ Hence
$\mathcal D k_D(a;e_1)\geq C$ and therefore, $L_{\mathcal
Dk_D}(\gamma)\ge C/2>0,$ which completes the proof of this
example.

Note that it shows that, with respect to the lengths of curves,
$\mathcal D k_D$ behaves different than the "real" derivative of
$k_D$ (cf. \cite{Ven} or \cite{Jar-Pfl2}, page 12). Moreover, it
implies that, in general, $\mathcal D k_D\neq\underline{\mathcal
D}k_D$.

\smallskip

\noindent{\bf Questions.} It will be interesting to know examples
showing that, in general, $\kappa_D\neq\mathcal D\tilde k_D$. It
remains also unclear whether $\mathcal D k_D$ is holomorphi\-cally
contractible (see \cite{Jar-Pfl}). Recall that $\int\mathcal D
k_D=k_D;$ but we do not know if $\int\underline{\mathcal D}k_D=k_D.$

\section{Proofs}

\noindent{\it Proof of Proposition \ref{pr2}.} First, we shall
consider the case $m=1$.  The key is the following

\begin{theorem}\label{th4}\footnote{We may replace Theorem
\ref{th4} by the approach used in the proof of the upper
semicontinuity of $\tilde k_M$.} \cite{Roy} Let $M$ be an
$n$-dimensional complex manifold and $f\in\mathcal O(\Bbb D,M)$
regular at $0$. Let $r\in(0,1)$ and $D_r=r\Bbb D\times\Bbb D^{n-1}$.
Then there exists $F\in\mathcal O(D_r,M)$, which is regular at $0$
and $F|_{r\Bbb D\times\{0\}}=f$.
\end{theorem}

Since $\kappa_M(z;0)=\mathcal D\tilde k_M(z;0)=0$, we may assume
that $X\neq 0$. Let $\alpha>0$ and $f\in\mathcal O(\Bbb D,M)$ be
such that $f(0)=z$ and $\alpha f_\ast(d/d\zeta)=X$. Let $r\in(0,1)$
and $F$ as in Theorem \ref{th4}. Since $F$ is regular at $0$, there
exist open neighborhoods $U=U(z)\subset M$ and $V=V(0)\subset D_r$
such that $F|_V:V\to U$ is biholomorphic. Hence $(U,\varphi)$ with
$\varphi=(F|_V)^{-1}$, is a chart near $z$. Note that
$\varphi_\ast(X)=\alpha e_1$, where $e_1=(1,0,\dots,0).$

If $w$ and $Y$ are sufficiently near $z$ and $\alpha e_1$,
respectively, then
$$
g(\zeta):=F(\varphi(w)+\zeta Y/\alpha),\quad \zeta\in r^2\Bbb D,
$$
belongs to $\mathcal O(r^2\Bbb D,M)$ with $g(0)=w$ and
$g(t\alpha)=\varphi^{-1}(\varphi(w)+tY),$ $t<r^2/\alpha.$ Therefore,
$r^2\widetilde k^\ast_M(w,\varphi^{-1}(\varphi(w)+tY))\leq t\alpha$.
Hence $r^2\mathcal D\tilde k_M(z;X)\le\alpha$.  Letting $r\to 1$ and
$\alpha\to\kappa_M(z;X)$ we get that $\mathcal D\tilde
k_M(z;X)\le\kappa_M(z;X).$

Let now $m\in\Bbb N$. By definition, $\kappa_M^{(m)}(z;\cdot)$ is
the largest function with the following property:

For any $X=\sum_{j=1}^mX_j$ one has that
$\kappa_M^{(m)}(z;X)\le\sum_{j=1}^m\kappa_M(z;X_j)$.

To prove that  $\kappa_M^{(m)}\ge\mathcal D k_M^{(m)}$ it suffices
to check that $\mathcal D k^{(m)}_M(z;\cdot)$ has the same property.
Following the above notation and choosing $Y_j\to\varphi_\ast X_j$
with $\sum_{j=1}^m Y_j=Y$, we set $w_0=w$ and
$w_j=\varphi^{-1}(\varphi(w)+t\sum_{k=1}^{j}Y_j)$. Since
$$k^{(m)}_M(w,w_q)\le\sum_{j=1}^m\tilde k_M(w_{j-1},w_j),$$ it follows
by the case $m=1$ that $$\mathcal D
k^{(m)}_M(z;X)\le\sum_{j=1}^m\mathcal D
k_M(z;X_j)\le\sum_{j=1}^m\kappa_M(z;X_j).$$

Finally, let $m=\infty$ and $n=\dim M$. Since
$\hat\kappa_M=\kappa_M^{(2n-1)}$ and $k_M\le k_M^{(2n-1)}$, we get
that $\mathcal D k_M\le\hat\kappa_M$ using the case $m=2n-1$. \qed

\smallskip

\noindent{\it Proof of Theorem \ref{th1}.} We may assume that $X\neq
0$. In virtue of Proposition \ref{pr2}, we have to show that
$$\kappa_M^{(m)}(z;X)\le\underline{\mathcal D}k^{(m)}_M(z;X).$$
For simplicity we assume that $M$ is a domain in $\Bbb C^n$.

(i) Fix a neighborhood $U=U(z)\Subset M.$ Applying the hyperbolicity
of $M$ at $z$, there are a neighborhood $V=V(z)\subset U$ and a
$\delta\in (0,1)$ such that, if $h\in\mathcal O(\Bbb D,M)$ with
$h(0)\in V$, then $h(\delta\Bbb D)\subset U$. Hence, by the Cauchy
inequalities, $||h^{(k)}(0)||\leq c/\delta^k$, $k\in\Bbb N.$

Now choose sequences $M\ni w_j\to z$, $\Bbb C_\ast\ni t_j\to 0$, and
$\Bbb C^n\ni Y_j\to X$ such that
$$
\frac{\widetilde k_M(w_j,w_j+t_jY_j)}{|t_j|}\to\underline{\mathcal
D}\widetilde k_M(z;X).
$$
There are holomorphic discs $g_j\in\mathcal O(\Bbb D,M)$ and
$\beta_j\in (0,1)$ with $g_j(0)=w_j$, $g_j(\beta_j)=w_j+t_jY_j,$ and
$\beta_j\le\widetilde k^\ast_M(w_j,w_j+t_jY_j)+|t_j|/j.$ Note that
$\widetilde k^\ast_M(w_j,w_j+t_jY_j)\le c_1||t_jY_j||\le c_2 |t_j|.$

Write
$$
w_j+t_jY_j=g_j(\beta_j)=w_j+g_j'(0)\beta_j+ h_j(\beta_j).
$$
Then $$||h_j(\beta_j)||\leq c\sum_{k=2}^\infty
\left(\tfrac{\beta_j}{\delta}\right)^k\leq c_3|\beta_j|^2\le
c_4|t_j|^2,\quad j\ge j_0.$$

Put $\widehat Y_j=Y_j-h_j(\beta_j)/t_j.$ We have that $g_j(0)=w_j$
and $\beta_jg_j'(0)/t_j=\widehat Y_j\to X$. Therefore,
$$ \kappa_M(w_j;\widehat
Y_j)\leq\frac{\beta_j}{|t_j|}\leq\frac{\widetilde
k^\ast_M(z_j,w_j+t_jY_j)}{|t_j|}+\frac{1}{j}.
$$
Hence with $j\to\infty$, we get that
$\kappa_M(z;X)=\underline{\kappa}_M(z;X)\leq \underline{\mathcal
D}\widetilde k_M(z;X)$.

(ii) The proof of the case $m\in\Bbb N$ is similar to the next one
and we omit it. Now, we shall consider the case $m=\infty$.

Note first that our assumption implies that $M$ is hyperbolic at $z$
and, by the contrary,
$$\forall\varepsilon>0\
\exists\delta>0:||w-z||<\delta,||Y-X||<\delta||X||$$
$$\Rightarrow|\kappa_M(w;Y)-\kappa_M(z;X)|<\varepsilon\kappa_M(z;X).\leqno{(1)}$$
Moreover, the proof of (i) shows that
$$\tilde k_M(a,b)\ge\kappa_M(a;b-a+o(a,b)),\hbox{ where }\lim_{a,b\to
z}\frac{o(a,b)}{||a-b||}=0.\leqno{(2)}$$

Choose now sequences $M\ni w_j\to z$, $\Bbb C_\ast\ni t_j\to 0$, and
$\Bbb C^n\ni Y_j\to X$ such that
$$\frac{k_M(w_j,w_j+t_jY_j)}{|t_j|}\to\underline{\mathcal
D}k_M(z;X).$$ There are points
$w_{j,0}=w_j,\dots,w_{j,m_j}=w_j+t_jX_j$ in $M$ such that
$$\sum_{k=1}^{m_j}\tilde k_M(w_{j,k-1},w_{j,k})\le
k_M(w_j,w_j+t_jY_j)+\frac{1}{j}.\leqno{(3)}$$ Set $w_{j,k}=w_j$ for
$k>m_j$. Since
$$k_M(w_j,w_{j,l})\le\sum_{j=1}^l\tilde k_M(w_{j,k-1},w_{j,k})\le
k_M(w_j,w_j+t_jY_j)+\frac{1}{j}\le c_2|t_j|+\frac{1}{j},$$ then
$k_M(w_j,w_{j,l})\to 0$ uniformly in $l.$ Then the hyperbolicity of
$M$ at $z$ implies that $w_{j,l}\to z$ uniformly in $l.$ Indeed,
assuming the contrary and passing to a subsequence, we may suppose
that $w_{j,l_j}\not\in U$ for some $U=U(z).$ Then
$$0=\lim_{j\to\infty}k_M(w_j,w_{j,l})\ge\ds\liminf_{z'\to z,w\in
M\setminus U}\tilde k_M(z',w)>0,$$ a contradiction.

Fix now $R>1.$ Then (1) implies that
$$\kappa_M(z;w_{j,k}-w_{j,k-1})\le R\kappa_M(w_{j,k};w_{j,k}-w_{j,k-1}
+o(w_{j,k},w_{j,k-1})),\ j\ge j(R).$$ It follows by this inequality,
(2) and (3) that
$$\sum_{k=1}^{m_j}\kappa_M(z;w_{j,k}-w_{j,k-1})
\le Rk_M(w_j,w_j+t_jYj)+\frac{R}{j}.$$ Since
$\hat\kappa_M(z;t_jY_j)$ is bounded by the first sum, we obtain that
$$\hat\kappa_M(z;Y_j)\le R\frac{k_M(w_j,w_j+t_jYj)+1/j}{|t_j|}.$$
Note that $\hat\kappa_M(z;\cdot)$ is a continuous function. Hence
with $j\to\infty$ and $R\to 1,$ we get that
$\hat\kappa_M(z;X)\le\underline{\mathcal D}k_M(z;X)$.\qed

\smallskip

\noindent{\bf Remark.} It follows by the above proofs and a standard
diagonal process that
$\underline{\kappa}_M(z;\cdot)=\underline{\mathcal D}\tilde
k(z;\cdot)$ if $M$ is hyperbolic at $z$.

\end{document}